\documentclass[12pt]{article}
\usepackage{amsthm,amsfonts,amssymb,amscd}

\begin{document}


\begin{center}
\large \bf Automorphisms of certain affine complements \\
in the projective space
\end{center}\vspace{0.5cm}

\centerline{A.V.Pukhlikov}\vspace{0.5cm}

\parshape=1
3cm 10cm \noindent {\small \quad\quad\quad \quad\quad\quad\quad
\quad\quad\quad {\bf }\newline We prove that every biregular
automorphism of the affine algebraic variety ${\mathbb
P}^M\setminus S$, $M\geqslant 3$, where $S\subset {\mathbb P}^M$
is a hypersurface of degree $m\geqslant M+1$ with a unique
singular point of multiplicity $(m-1)$, resolved by one blow up,
is a restriction of some automorphism of the projective space
${\mathbb P}^M$, preserving the hypersurface $S$; in particular,
for a general hypersurface $S$ the group $\mathop{\rm Aut}
({\mathbb P}^M\setminus S)$ is trivial.

Bibliography: 24 titles.} \vspace{1cm}

\vspace{1cm}

{\bf 1. Statement of the main result.} Let ${\mathbb P}= {\mathbb
P}^M$ be the complex projective space of dimension $M\geqslant 3$
and $S\subset{\mathbb P}$ a hypersurface of degree $m\geqslant
M+1$ with a unique singular point $o\in S$ of multiplicity $m-1$,
that can be resolved by one blow up. More precisely, let
$\sigma\colon{\mathbb P}^+\to{\mathbb P}$ be the blow up of the
point $o$ with the exceptional divisor
$E=\sigma^{-1}(o)\cong{\mathbb P}^{M-1}$. We assume that the
strict transform $S^+\subset{\mathbb P}^+$ is a non-singular
hypersurface and the projectivised tangent cone $S^+\cap E$ is a
non-singular hypersurface of degree $m-1$ in $E\cong{\mathbb
P}^{M-1}$. The main result of the present paper is the following
claim.\vspace{0.1cm}

{\bf Theorem 1.} {\it Every automorphism $\chi$ of the affine
algebraic set ${\mathbb P}\backslash S$ is the restriction of some
projective automorphism $\chi_{{\mathbb P}}\in\mathop{\rm
Aut}{\mathbb P}$, preserving the hypersurface $S$. In particular,
the group
$$
\mathop{\rm Aut}({\mathbb P}\backslash S)
$$
is finite and trivial for a Zariski general hypersurface
$S$.}\vspace{0.1cm}

Due to certain well known facts about automorphisms of projective
hypersurfaces (see, for instance, \cite{MatMon}) Theorem 1 is
easily implied (see Sec. 6) by a somewhat more general fact. Let
$S'\subset{\mathbb P}$ be one more hypersurface of degree $m$ with
a unique singular point $o'\in S'$ of multiplicity $m-1$, that is
resolved by one blow up (in the sense specified above). Then the
following claim is true.\vspace{0.1cm}

{\bf Theorem 2.} {\it Every isomorphism of affine algebraic
varieties
$$
\chi\colon{\mathbb P}\backslash S\to{\mathbb P}\backslash S'
$$
is the restriction of some projective automorphism $\chi_{{\mathbb
P}}\in\mathop{\rm Aut}{\mathbb P}$, transforming the hypersurface
$S$ into hypersurface $S'$.}\vspace{0.1cm}

Obviously, $\chi_{\mathbb P} (o)=o'$. It is Theorem 2 that we
prove below.\vspace{0.1cm}

If $z_1,\dots,z_M$ is a system of affine coordinates on ${\mathbb
P}$ with the origin at the point $o$, then the hypersurface $S$ is
defined by the equation
\begin{equation}\label{8}
f(z_*)=q_{m-1}(z_*)+q_m(z_*)=0,
\end{equation}
where $q_i(z_*)$ are homogeneous polynomials of degree $i$ in the
coordinates $z_*$. An irreducible hypersurface of that type is
rational and it is this property that makes the problem of
describing the group of automorphisms $\mathop{\rm Aut}({\mathbb
P}\backslash S)$ meaningful, see the discussion in Subsection 2
below.\vspace{0.1cm}

The paper is organized in the following way. In Sec. 2 we discuss
the general problem of describing the automorphisms of affine
complements and what little is known in that direction (for
non-trivial cases), and also some well known conjectures and
non-completed projects. In Sec. 3 we start the proof of Theorem 2:
for an arbitrary isomorphism of affine varieties
$$
\chi\colon{\mathbb P}\backslash S\to{\mathbb P}\backslash S'
$$
we define the key numerical characteristics (such as the ``degree)
and obtain the standard relations between them (for instance, an
analog of the ``Noether-Fano inequality'' for the affine case). In
Sec. 4 we construct the resolution of the maximal singularity of
the map $\chi$, which is now considered a birational map (Cremona
transformation) $\chi_{{\mathbb P}}\colon {\mathbb P}
\dashrightarrow{\mathbb P}$, the restriction of which onto the
affine complement ${\mathbb P}\backslash S$ is an isomorphism onto
${\mathbb P}\backslash S'$. Finally, in Sec. 5 we exclude the
maximal singularity, which completes the proof of Theorem
2.\vspace{0.1cm}

The author thanks the referee for useful comments and suggestions.
\vspace{0.3cm}


{\bf 2. Automorphisms of affine complements.} Let $X$ be a
non-singular projective rationally connected variety, $Y$ and $Y'$
irreducible ample divisors, so that their complements $X\backslash
Y$ and $X\backslash Y'$ are affine varieties. Two natural
questions can be asked:\vspace{0.1cm}

1) are the affine varieties $X\backslash Y$ and $X\backslash Y'$
isomorphic,\vspace{0.1cm}

2) if $Y'=Y$, what is the group of biregular automorphisms
$\mathop{\rm Aut}(X\backslash Y)$.\vspace{0.1cm}

It is natural to consider a biregular isomorphism $\chi\colon
X\backslash Y\to X\backslash Y'$ (if they exist) as a {\it
birational automorphism} $\chi_X\in\mathop{\rm Bir}X$, regular on
the affine open set $X\backslash Y$ and mapping it isomorphically
onto $X\backslash Y'$. The case when $\chi_X\in\mathop{\rm Aut}X$
is a biregular automorphism of the variety $X$ and the
corresponding isomorphism $\chi$ of affine complements itself we
will say to be trivial. We therefore consider the following
problem: are there any {\it non-trivial isomorphisms} $\chi\colon
X\backslash Y\to X\backslash Y'$, when
$$
\chi_X\in\mathop{\rm Bir}X\backslash\mathop{\rm Aut}X,
$$
and, respectively, are the groups
$$
\mathop{\rm Aut}(X\backslash Y)\,\,\mbox{and}\,\,\mathop{\rm
Aut}(X)_Y
$$
the same (the second symbol means the stabilizer of the divisor
$Y$ in the group $\mathop{\rm Aut}(X)$). In particular, Theorem 1
claims that
$$
\mathop{\rm Aut}({\mathbb P}\backslash S)=\mathop{\rm
Aut}({\mathbb P})_S
$$
for hypersurfaces $S\subset{\mathbb P}$, described in Sec.
1.\vspace{0.1cm}

{\bf Proposition 1.} {\it Let $\chi$ be a non-trivial isomorphism
of affine complements $X\backslash Y$ and $X\backslash Y'$. Then
$Y$ and $Y'$ are birationally ruled varieties, that is to say, for
some irreducible varieties $Z$ and $Z'$ of dimension $\mathop{\rm
dim}X-2$ the varieties $Y$ and $Y'$ are birational to the direct
products $Z\times{\mathbb P}^1$ and $Z'\times{\mathbb P}^1$,
respectively.}\vspace{0.1cm}

{\bf Proof.} The birational map $\chi^{-1}_X$ is regular at the
generic point of the divisor $Y'$, and its image can not be the
generic point of the divisor $Y$: in such case $\chi_X$ would have
been an isomorphism in codimension 1 and for that reason a
biregular automorphism, contrary to our assumption. Therefore,
$(\chi^{-1}_X)_*Y'\subset X$ is an irreducible subvariety of
codimension at least 2 (which is, of course, contained in $Y$).
Now let us consider a resolution of singularities
$\varphi\colon\widetilde{X}\to X$ of the map $\chi_X$. By what we
said above, there is an exceptional divisor
$E\subset\widetilde{X}$ of this resolution, such that
$$
(\chi\circ\varphi)|_E\colon E\to Y'
$$
is a birational map. Therefore, $Y'$ is a birationally ruled
variety. For $Y$ we argue in a symmetric way. Q.E.D. for the
proposition.\vspace{0.1cm}

{\bf Remark 1.} Assume in addition that $Y'$ is a rationally
connected variety. Then in the notations of the proof of
Proposition 1 we conclude that the centre of the exceptional
divisor $E$ on $X$, that is, the irreducible subvariety
$\varphi(E)$, is a rationally connected variety. This is also true
for the centres of the divisor $E$ on the ``lower'' stroreys of
the resolution of $\varphi$.\vspace{0.1cm}

{\bf Example 1.} There are no non-trivial isomorphisms of affine
complements ${\mathbb P}\backslash Y$ and ${\mathbb P}\backslash
Y'$, if $Y\subset{\mathbb P}$ is a non-singular hypersurface of
degree at least $M+1$. Indeed, the hypersurface $Y$ is not a
birationally ruled variety. \vspace{0.1cm}

Assume now that $Y\subset X$ is a Fano variety. It is well known
(see, for instance, \cite[Chapter 2]{Pukh13a}), that birationally
rigid Fano varieties are not birationally ruled. Therefore, if $Y$
is a birationally rigid variety, then every isomorphism of affine
complements $X\backslash Y$ and $X\backslash Y'$ (where $Y'\subset
X$ is an irreducible ample divisor) is trivial, that is, it
extends to an automorphism of the variety $X$. In particular, the
group $\mathop{\rm Aut}(X\backslash Y)$ is $\mathop{\rm
Aut}(X)_Y$. This makes it possible to construct numerous examples
of affine complements with no non-trivial isomorphisms and
automorphisms. Below we give some of them.\vspace{0.1cm}

{\bf Example 2.} Let $V\subset{\mathbb P}^4$ be a smooth
three-dimensional quartic. Because of its birational superrigidity
(\cite{IM}) the affine complement ${\mathbb P}^4\backslash V$ has
no non-trivial automorphisms and isomorphisms. The same is true
for quartics with at most isolated double points, provided that
the variety $V$ is factorial and its singularities are terminal,
see \cite{Pukh89c,Me04} and \cite{Shr08b}.\vspace{0.1cm}

{\bf Example 3.} Let $V\subset{\mathbb P}^M$ be a general smooth
hypersurface of degree $M$, where $M\geqslant 5$. Because of its
birational superrigidity (\cite{Pukh98b}) the affine complement
${\mathbb P}\backslash V$ has only trivial automorphisms and
isomorphisms. The same is true if we allow $V$ to have quadratic
singularities of rank at least 5 \cite{EP}. This example
generalizes naturally for Fano complete intersections. Let
$k\geqslant 2$,
$$
Y=F_1\cap\dots\cap F_k\subset{\mathbb P}^{M+k}
$$
be a non-singular complete intersection of codimension $k$, where
$F_i$ is a hypersurface of degree $d_i$, and
$$
d_1+\dots+d_k=M+k,
$$
$M\geqslant 4$, that is, $Y$ is a non-singular $M$-dimensional
Fano variety of index 1. For $i\in\{1,\dots,k\}$ set
$$
X_i=\bigcap_{j\neq i}F_j
$$
and assume that $X_i$ is also non-singular. Then
$X_i\subset{\mathbb P}^{M+k}$ is a $(M+1)$-dimensional Fano
variety of index $d_i+1$, containing $Y$ as a very ample divisor,
so that the complement $X_i\backslash Y$ is an affine variety. If
the set of integers $(d_1,\dots,d_k)$ satisfies the conditions of
any of the papers \cite{Pukh01,Pukh14a,Pukh13c}, and the variety
$Y$ is sufficiently general in its family, then due to its
birational superrigidity the equality
$$
\mathop{\rm Aut}(X_i\backslash Y)=\mathop{\rm Aut}(X_i)_Y
$$
holds (and a similar claim for automorphisms). Of course, these
arguments are non-trivial only for those cases, when $\mathop{\rm
Aut}X_i\neq\mathop{\rm Bir}X_i$: for instance, for $k=2$ and
$(d_1,d_2)=(2,M)$ the variety $X_2$ is a $(M+1)$-dimensional
quadric and its group of birational automorphisms is the Cremona
group of rank $M+1$. We get another non-trivial example for $k=2$
and $(d_1,d_2)=(3,M-1)$, where $X_2$ is a $(M+1)$-dimensional
cubic hypersurface which has a huge group of birational
automorphisms. Using other families of birationally superrigid or
rigid Fano varieties, one can construct more non-trivial examples
of affine complements, all automorphisms of which are
trivial.\vspace{0.1cm}

{\bf Example 4.} In \cite{Kol95a} it is shown that a very general
hypersurface $V_d\subset{\mathbb P}$ for $d\geqslant\frac23M$ is
not birationally ruled. Therefore, for such hypersurfaces their
affine complements ${\mathbb P}\backslash V_d$ have no non-trivial
isomorphisms and automorphisms.\vspace{0.1cm}

{\bf Example 5.} In \cite{Pukh15a} it was shown that a Zariski
general hypersurface $V_{M-1}\subset{\mathbb P}$ for $M\geqslant
16$ has no other structures of a rationally connected fibre space
apart from pencils of hyperplane sections. In particular, it has
no structures of a conic bundle and for that reason is not
birationally ruled. It follows that for for those hypersurfaces
the affine complements ${\mathbb P}\backslash V_{M-1}$ have no
non-trivial isomorphisms and automorphisms.\vspace{0.1cm}

Unfortunately, if the variety $Y$ is birationally ruled, then the
problem of describing the isomorphisms of the affine complement
$X\backslash Y$ becomes very hard (except for trivial cases, when,
for instance, the variety $X$ itself satisfies the equality
$\mathop{\rm Bir}X=\mathop{\rm Aut}X)$. The only complete result
here is Theorem 2 of the present paper. As for the main objects of
study today, they are particular classes of three-dimensional
affine complements, such as the complement ${\mathbb
P}^3\backslash S$ to a cubic surface (non-singular or with
prescribed singularities) or the affine space ${\mathbb A}^3$ and
certain similar affine varieties. In respect of complements to
cubic surfaces there is a classical conjecture, stated by
M.Kh.Gizatullin in \cite[p. 6]{Giz2005}: if the cubic surface $S$
is non-singular, then its complement ${\mathbb P}^3\backslash S$
has no non-trivial automorphisms. However, if the cubic surface
has a double point, then non-trivial automorphisms do exist ---
they were discovered by S.Lamy and J.Blanc (as far as the author
knows, those examples were not published). A similar conjecture
was stated by A.Dubulouz for the case when $X$ is a Fano double
cover of index 2, branched over a surface $W\subset{\mathbb P}^3$
of degree 4 and  $S$ is the inverse image of a plane in ${\mathbb
P}^3$.\vspace{0.1cm}

In respect of the groups of automorphisms of affine varieties a
huge material has been accumulated; there are a lot of results
about special groups of automorphisms, dynamical properties of
particular automorphisms etc. We only point out three recent
papers \cite{Giz2008,Lamy2013,Lamy2014}, see also the bibliography
in those papers.\vspace{0.1cm}

The groups of automorphisms of affine algebraic surfaces are much
better understood: here we have such fundamental results as the
complete description of the groups of automorphisms of the plane
$\mathop{\rm Aut}{\mathbb A}^2$, see \cite{GizDan1975,GizDan1977}.
This direction is still being actively explored
\cite{DubLamy2015,FurLamy2010,CanLamy2006,Lamy2005}.\vspace{0.3cm}


{\bf 3. Start of the proof of Theorem 2.} Let
$$
\chi\colon{\mathbb P}\backslash S\to{\mathbb P}\backslash S'
$$
be an isomorphism of affine varieties. Assume that $\chi$ is
non-trivial, that is, the corresponding birational map
$\chi_{\mathbb P}\colon{\mathbb P}\dashrightarrow{\mathbb P}$ is
not a biregular isomorphism. Let
$$
\varphi\colon\widetilde{\mathbb P} \to {\mathbb P}
$$
be its resolution (a sequence of blow ups with non-singular
centres), so that $\psi=\chi_{\mathbb P}
\circ\varphi\colon\widetilde{\mathbb P} \to{\mathbb P}$ is
aregular map. Furthermore, set ${\cal E}_{\varphi}$ to be the set
of prime $\varphi$-exceptional divisors. By assumption, for the
strict transform of $S'$ we have
$$
T=(\chi_{\mathbb P} \circ\varphi)^{-1}_*S'\in{\cal E}_{\varphi}.
$$
Ser $B=\varphi(T)$ to be the centre of the exceptional divisor $T$
on ${\mathbb P}$, an irreducible subvariety of codimension at
least 2, and moreover $B\subset S$. Therefore we get the positive
integers $a=a(T,{\mathbb P})$ (the discrepancy of the divisor $T$
with respect to ${\mathbb P}$) and
$$
b=\mathop{\rm ord}\nolimits_TS=\mathop{\rm
ord}\nolimits_T\varphi^*S.
$$
Furthermore, let $\Sigma$ be the strict transform of the linear
system of hyperplanes with respect to $\chi_{\mathbb P}$. This is
a mobile linear system $\Sigma\subset|nH|$, where $H$ is a
hyperplane in ${\mathbb P}$, and $n\geqslant 2$. Set
$$
\nu=\mathop{\rm ord}\nolimits_T\varphi^*\Sigma.
$$

{\bf Proposition 2.} {\it The following equalities are true:}
\begin{equation}\label{1}
bn=\nu m+1,\quad (M+1)b=am+(M+1).
\end{equation}

{\bf Proof.} Write down ${\cal E}={\cal E}_{\varphi}\setminus
\{T\}$, so that ${\cal E}_{\varphi}={\cal E}\coprod\{T\}$. Let
$D\in\Sigma$ be a general divisor,
$\widetilde{D}\in\widetilde{\Sigma}$ its strict transform on
$\widetilde{\mathbb P}$, where $\widetilde{\Sigma}$ is the strict
transform of the linear system $\Sigma$ on $\widetilde{\mathbb P}$
with respect to $\varphi$. Let
$\widetilde{S}\subset\widetilde{\mathbb P}$ be the strict
transform of the hypersurface $S$. By the symbol $\widetilde{K}$
we denote the canonical class of the variety $\widetilde{\mathbb
P}$. We obtain the following presentations:
$$
\begin{array}{ccccccc}
\widetilde{D} &\sim & nH & - & \nu T & - & \sum\limits_{E\in{\cal E}}\nu_EE,\\
  &  &  &  &  &  & \\
\widetilde{K} & = & -(M+1)H & + & aT & + & \sum\limits_{E\in{\cal E}}a_EE,\\
  &  &  &  &  &  & \\
\widetilde{S} & \sim & mH & - & bT & - & \sum\limits_{E\in{\cal E}}b_EE,\\
\end{array}
$$
where the coefficients $\nu_E, a_E, b_E$ have the obvious meaning
(in order to simplify the formulas we write $H$ in stead of
$\varphi^*H$). Consider the family of lines ${\cal L}$ on
${\mathbb P}$. Obviously, a general line $L\in {\cal L}$ does not
meet the set
$$
\mathop{\bigcup}\limits_{E\in{\cal
E}}\psi(E)\cup\psi(\widetilde{S}),
$$
since it is of codimension at least 2 (recall that
$\widetilde{S}\subset\widetilde{\mathbb P}$ is a
$\psi$-exceptional divisor). Therefore, the strict transform
$\widetilde{L}\subset\widetilde{\mathbb P}$ satisfies the
equalities
\begin{equation}\label{2}
(\widetilde{L}\cdot\widetilde{D})=1,\,\,(\widetilde{L}\cdot\widetilde{K})=
-(M+1),\,\,(\widetilde{L}\cdot\widetilde{S})=0.
\end{equation}
Besides, $(\widetilde{L}\cdot T)=(L\cdot S')=m$. Set
$$
d=(\widetilde{L}\cdot H).
$$
Obviously, $d$ is the degree of the curve
$\varphi(\widetilde{L})\subset\widetilde{\mathbb P}$ in the usual
sense. Finally, we have the equality $(\widetilde{L}\cdot E)=0$
for every exceptional divisor $E\in{\cal E}$. Therefore the
equalities (\ref{2}) imply the relations
$$
dn-\nu m=1,
$$
$$
-(M+1)d+am=-(M+1),
$$
$$
dm-bm=0.
$$
The last equality implies that $d=b$. Now the equalities (\ref{1})
follow in a straightforward way. Q.E.D. for the
proposition.\vspace{0.1cm}

{\bf Remark 2.} The relations (\ref{1}) imply the equality
$$
\nu=\frac{a}{M+1}n+\frac{n-1}{m}.
$$
Since $n\geqslant 2$, we obtain the inequality
$$
\nu>\frac{a}{M+1}n.
$$
This is the usual Noether-Fano inequality for the birational map
$\chi_{\mathbb P}$. Therefore, the prime divisor $T$ (the strict
transform of the hypersurface $S'$ on $\widetilde{\mathbb P}$) is
a {\it maximal singularity} of the linear system $\Sigma$ (see,
for instance, Definition 1.4 in \cite[Chapter
2]{Pukh13a}).\vspace{0.1cm}

Although the relations (\ref{1}) are sufficient for the proof of
Theorem 2, we will show similar relations for every infinitely
near divisor $E\in{\cal E}$. Recall that we defined the integers
$$
a_E=a(E,{\mathbb P})\quad \mbox{è}\quad b_E=\mathop{\rm
ord}\nolimits_E\varphi^*S,
$$
where the discrepancy is understood with respect to the birational
morphism $\varphi$. Let $a'_E$ be the discrepancy of the divisor
$E$ with respect to the birational morphism $\psi$ and
$b'_E=\mathop{\rm ord}_E\psi^*S'$, so that we get the equality
\begin{equation} \label{5}
\widetilde{K}=\psi^*K_{\mathbb P}+a'\widetilde{S}+\sum_{E\in{\cal
E }}a'_EE
\end{equation}
and the presentation
$$
\widetilde{S}'=T\sim\psi^*(mH)-b'\widetilde{S}-\sum_{E\in{\cal
E}}b'_EE,
$$
where $a'>0$ and $b'>0$ have the same sense in respect of the
image of the map $\chi$ as $a$ and $b$ in respect of the original
projective space ${\mathbb P}$.\vspace{0.1cm}

{\bf Proposition 3.} {\it For every divisor $E\in{\cal E}$ the
following equalities hold:}
\begin{equation}\label{3}
(M+1)b_E+ma'_E=(M+1)b'_E+ma_E
\end{equation}
{\it and}
\begin{equation}\label{4}
b_En=m\nu_E+b'_E.
\end{equation}

{\bf Proof.} Consider a mobile family of curves ${\cal C}$ on
${\mathbb P}$ with the following properties:\vspace{0.1cm}

1) every curve $C\in{\cal C}$ is an irreducible rational curve of
degree $l\geqslant 2$,\vspace{0.1cm}

2) the strict transform $\widetilde{C}$ of a general curve
$C\in{\cal C}$ on $\widetilde{\mathbb P}$ with respect to the
birational morphism $\psi$ meets $E$ transversally at a unique
point $p_C$ of general position on $E$ and does not meet other
$\psi$-exceptional divisors, in particular
$\widetilde{C}\cap\widetilde{S}=\emptyset$,\vspace{0.1cm}

3) the curves of the family ${\cal C}$ sweep out a Zariski open
subset of the space ${\mathbb P}$.\vspace{0.1cm}

Such a family of rational curves is easy to construct using the
methods of elementary algebraic geometry, see \cite[Chapter 2,
Section 3]{Pukh13a}. Let $p\in E$ be a point of general position,
$q=\psi(p)\in{\mathbb P}$ its image on ${\mathbb P}$ and
$(v_1,\dots,v_M)$ a system of affine coordinates on ${\mathbb P}$
with the origian at that point. We construct the curve $C$ in the
parametric form:
$$
\begin{array}{lcccccc}
v_1 & = & \alpha_{1,1}t & + & \dots & + & \alpha_{1,l}t^l,\\
& & & & \dots & & \\
v_M & = & \alpha_{M,1}t & + & \dots & + & \alpha_{M,l}t^l,\\
\end{array}
$$
where $l$ is sufficiently large. In \cite{Pukh95a}, see also
\cite[Chapter 2, Theorem 3.1]{Pukh13a} it was shown that there is
a set of coefficients $\alpha_{i,j},i=1,\dots,M$, $j=1,\dots,a'_E$
(in fact, in stead of $a'_E$ one can take an essentially smaller
number, but we do not need that), such that for any coefficients
$\alpha_{i,j}$ for $j\geqslant a'_E+1$ the strict transform of
such curve meets $E$ transversally at the point $p$ when $t=0$.
Varying the coefficients $\alpha_{i,j}$ for $j\geqslant a'_E+1$,
one can ensure that the curve $C$ goes through the point $q$ only
when $t=0$ and intersects the closed subset of codimension
$\geqslant 2$
$$
\bigcup\limits_{E\in{\cal E}}\psi(E)\cup\psi(\widetilde{S})
$$
only at the point $q$. Such curves satisfy the properties
1)-3).\vspace{0.1cm}

Now we argue in exactly the same way as in the proof of
Proposition 2. We have the equality
$$
(\widetilde{C}\cdot\widetilde{D})=l=dn-\nu(lm-b'_E)-\nu_E,
$$
where $d=\mathop{\rm deg}\varphi(\widetilde{C})$. Multiplying
$\widetilde{C}$ by the canonical class $\widetilde{K}$ and using
the presentation (\ref{5}), we obtain the equality
$$
-(M+1)d+a(lm-b'_E)+a_E=-(M+1)l+a'_E.
$$
Finally, multiplying the curve $\widetilde{C}$ by $\widetilde{S}$,
we get:
$$
dm-b(lm-b'_E)-b_E=0
$$
(the expression in brackets $lm-b'_E$ is the ``residual
intersection'' $(\widetilde{C}\cdot T$)). From here, using the
equalities (\ref{1}), by means of easy computations we get the
equalities (\ref{3}) and (\ref{4}). Q.E.D. for the
proposition.\vspace{0.3cm}


{\bf 4. The resolution of the maximal singularity.} Let
$$
\varphi_{i,i-1}\colon X_i\to X_{i-1},
$$
$i=1,\dots,N$, be the resolution of the maximal singularity $T$ of
the linear system $\Sigma$ in the sense of \cite[Chapter
2]{Pukh13a}, that is, $X_0={\mathbb P}$, each map
$\varphi_{i,i-1}$ is the blow up of the (possibly singular)
irreducible subvariety $B_{i-1}\subset X_{i-1}$, which is the
centre of the exceptional divisor $T$ on $X_{i-1}$. Set
$E_i=\varphi^{-1}_{i,i-1}(B_{i-1})$. The strict transform of the
subvariety $R\subset X_i$ on a higher storey of the resolution
$X_j$, where $j>i$, is denoted by adding the upper index $j$: we
write $R^j$.\vspace{0.1cm}

For $j>i$ we set
$$
\varphi_{j,i}=\varphi_{i+1,i}\circ\dots\circ\varphi_{j,j-1}\colon
X_j\to X_i.
$$
The exceptional divisor $E_N\subset X_N$ of the last blow up
realizes the maximal singularity $T$: the birational map
$$
\varphi^{-1}\circ\varphi_{N,0}\colon
X_N\dashrightarrow\widetilde{\mathbb P}
$$
is regular at the general point of the divisor $E_N$ and maps it
onto $T$.\vspace{0.1cm}

On the set $\{1,\dots,N\}$ there is a natural structure of an
oriented graph: $i\to j$, if and only if $i>j$ and the inclusion
$$
B_{i-1}\subset E^{i-1}_j
$$
holds. If the vertices $i$ and $j$ are not joined by an orinted
edge, we write $i\nrightarrow j$.\vspace{0.1cm}

For $i\neq j$ we denote by the symbol $p_{ij}$ the number of paths
in that graph from the vertex $i$ to the vertex $j$ (so that
$p_{ij}=0$ for $i<j$ and $p_{ij}\geqslant 1$ for $i>j$). For
convenience we set $p_{ii}=1$ for $i=1,\dots,N$. Finally, in order
to simplify our notations, we write $p_i$ in stead of $p_{Ni}$.
Let
$$
\delta_i=\mathop{\rm codim}B_{i-1}-1
$$
be the elementary discrepancies. Then the following equality
holds:
$$
a=\sum^N_{i=1}p_i\delta_i.
$$
Let us also introduce the elementary multiplicities
$$
\nu_i=\mathop{\rm mult}\nolimits_{B_{i-1}}\Sigma^{i-1},
$$
$i=1,\dots,N$ (where, in accordance with the general principle of
notations, $\Sigma^{i-1}$ means the strict transform of the mobile
linear system $\Sigma$ on $X_{i-1}$) and
$$
\mu_i=\mathop{\rm mult}\nolimits_{B_{i-1}}S^{i-1},
$$
$i=1,\dots,N$. Obviously,
$$
\nu=\sum^N_{i=1}p_i\nu_i\quad \mbox{and} \quad
b=\sum^N_{i=1}p_i\mu_i.
$$
Note that for some $k\leqslant N$ the strict transform $S^{k-1}$
contains $B_{k-1}$, but $S^k$ no longer contains $B_k$, so that
$\mu_{k+1}=\dots=\mu_N=0$, and for that reason
$$
b=\sum^k_{i=1}p_i\mu_i.
$$

If $B_0\neq o$ is not the unique singular point of the
hypersurface $S$, then obviously
$$
\mu_1=\dots=\mu_k=1,
$$
so that $b=p_1+\dots+p_k$. If $B_0=o$, then $\mu_1=m-1$ and by the
assumption about the singularities of the divisor $S$ the strict
transform $S^1$ is smooth, so that $\mu_i=1$ for $2\leqslant
i\leqslant k$. Therefore if $B_0=o$, then
$$
b=(m-1)p_1+p_2+\dots+p_k.
$$

Finally, let us point out one property of the numbers $p_i$. Since
by construction we have $\varphi_{i,i-1}(B_i)=B_{i-1}$ ($B_i$ is
the centre of the exceptional divisor $T$ on $X_i$, and $B_{i-1}$
is its centre on $X_{i-1}$), the dimensions $\mathop{\rm dim}B_i$
do not decrease when $i$ is growing. Accordingly, the codimensions
$\mathop{\rm codim}B_i$ do not increase when $i$ is growing, so
that $\delta_1\geqslant\delta_2\geqslant\dots\geqslant\delta_N$.
Assume that for some $k_1<k$ the centres of the blow ups
$B_{i-1}$, $i=k_1+1,\dots,k$, have the maximal dimension
$M-2$.\vspace{0.1cm}

{\bf Proposition 4.} {\it Under the assumptions above for
$k-k_1\geqslant 3$ the subgraph with the vertices $k_1+1,\dots,k$
is a chain:
$$
k_1+1\leftarrow k_1+2\leftarrow\dots\leftarrow
k,
$$
that is, between the vertices of the subgraph there are no other
arrows apart from the consecutive ones $i\leftarrow i+1$.
Moreover,}
$$
p_{k_1+1}=\dots=p_{k-1}.
$$

{\bf Proof.} By the definition of the number $k$, for $i\leqslant
k$ we have $B_{i-1}\subset S^{i-1} $, where the divisor $S^{i-1}$
is non-singular at the general point $B_{i-1}$ for $i\geqslant
k_1+1$. Therefore, for $k_1+1\leqslant i\leqslant k-2$ we have
$$
B_i=E_i\cap S^i\quad\mbox{and}\quad B_{i+1}=E_{i+1}\cap S^{i+1}
$$
(since $B_i$ is contained in both $E_i$ and $S^i$ and has
codimension 2, and the same is true for $B_{i+1}$), and $E_i$ and
$S^i$ (respectively, $E_{i+1}$ and $S^{i+1}$) meet transversally
at the general point of $B_i$ (respectively, of $B_{i+1}$), so
that $E^{i+1}_i$ and $S^{i+1}$ do not meet over a point of general
position in $B_i$. Therefore, $B_{i+1}\not\subset E^{i+1}_i$ and
the first claim of the proposition is shown.\vspace{0.1cm}

In particular, $k\nrightarrow k-2$. But then for any vertex
$e\geqslant k+1$ we have $e\nrightarrow k-2$, either, so that
every path from the vertex $N$ to the vertex $i\leqslant k-2$ must
go through the vertex $k-1$. This proves the second claim of
Proposition 4. Q.E.D.\vspace{0.1cm}

Now we are ready to complete the proof of Theorem 2.\vspace{0.3cm}


{\bf 5. Exclusion of the maximal singularity.} Let us write down
the second of the equalities (\ref{1}) in terms of elementary
multiplicities and discrepancies:
\begin{equation}\label{6}
(M+1)\sum^k_{i=1}p_i\mu_i=m\sum^N_{i=1}p_i\delta_i+(M+1).
\end{equation}
We conclude immediately that $B_0=o$ is the singular point of the
hypersurface $S$. Otherwise all multiplicities $\mu_i=1$, so that
from the formula (\ref{6}) we would have obtained
$$
0=\sum^k_{i=1}(m\delta_i-M-1)p_i+m\sum^N_{i=k+1}p_i\delta_i+(M+1),
$$
which is impossible, since $m\geqslant M+1$ and $\delta_i\geqslant
1$, so that all three components in the right hand side of the
last formula are non-negative and at least one of them is
positive.\vspace{0.1cm}

So $B_0=o$. Here $\mu_1=m-1$ and $\delta_1=M-1$, so we get the
equality
\begin{equation}\label{7}
(2m-M-1)p_i=
\sum^k_{i=2}(m\delta_i-M-1)p_i+m\sum^N_{i=k+1}p_i\delta_i+(M+1),
\end{equation}
all components in which both in the right and left hand side are
non-negative. By Remark 1, all centres $B_i$ of the blow ups
$\varphi_{i+1,i}$ are rationally connected varieties. In
particular, $B_1\neq S^1\cap E_1$, since $S^1\cap E_1\subset
E_1\cong{\mathbb P}^{M-1}$ is a non-singular hypersurface of
degree $m-1\geqslant M$, which is not rationally connected. Thus
if $k\geqslant 2$, then $B_1$ is a subvariety of codimension at
least 3 in $X_1$, so that $\delta_2\geqslant 2$ and the
coefficient at $p_2$ is not smaller than $2m-M-1$. If also
$3\nrightarrow 1$, then $p_1=p_2$ (every path from the vertex $N$
to the vertex 1 must go through the vertex 2) and we ob tain a
contradiction: in the equality (\ref{7}) the right hand side is
strictly higher than the left hand side. If $N=2$, then
$p_1=p_2=1$ and we obtain a contradiction again: in this case the
equality (\ref{7}) takes the form
$$
2m-M-1=\delta_2m-M-1+(M+1)
$$
with $\delta_2\geqslant 2$, which is also impossible. We conclude
that for $k\geqslant 2$ with necessity $N\geqslant 3$ and $3\to
1$.\vspace{0.1cm}

{\bf Proposition 5.} {\it The case $k=1$ is
impossible.}\vspace{0.1cm}

{\bf Proof.} Assume the converse: $k=1$. Then $b=(m-1)p_1$. Let
$Q\subset{\mathbb P}$ be a general hypersurface of degree $m$ with
the point $o$ as a singular point of multiplicity $m-1$. Since it
is general, $B_1\not\subset Q^1$, so that
$$
\mathop{\rm ord}\nolimits_T\varphi^*Q=\mathop{\rm
ord}\nolimits_{E_N}\varphi^*_{N,0}Q=(m-1)p_1=b,
$$
and it follows that for the strict transform
$\widetilde{Q}\subset\widetilde{\mathbb P}$ we get the
presentation
\begin{equation}\label{9}
\widetilde{Q}\sim mH-bT-\sum_{E\in{\cal E}}q_E E.
\end{equation}
Therefore, $(\widetilde{L}\cdot\widetilde{Q})=0$, where
$\widetilde{L}$ is the strict $\psi$-transform of a general line
$L\subset{\mathbb P}$ (see the proof of Proposition 2). But the
curves $\widetilde{L}$ sweep out a Zariski open subset of the
space ${\mathbb P}$, and the hypersurfaces $Q$ sweep out $\mathbb
P$. This contradiction proves Proposition 5.\vspace{0.1cm}

Set $l=\mathop{\rm max}\{i\,|\,i\to 1,1\leqslant i\leqslant
N\}$.\vspace{0.1cm}

{\bf Proposition 6.} {\it The case $l\leqslant k$ is
impossible.}\vspace{0.1cm}

{\bf Proof.} Assume the converse: $l\leqslant k$. We could see
above that $N\geqslant 3$ and $3\to 1$, so that $l\geqslant 3$.
For any $i\leqslant l$, $i\geqslant 2$, we have
$$
B_{i-1}\subset S^{i-1}\cap E^{i-1}_1\cap E_{i-1},
$$
so that $\mathop{\rm codim}B_{i-1}\geqslant 3$ and
$\delta_i\geqslant 2$. Let us re-write the right hand side of the
equality (\ref{7}) in the form
$$
\sum^l_{i=2}(m\delta_i-M-1)p_i+\sum^k_{i=l+1}
(m\delta_i-M-1)p_i+m\sum^N_{i=k+1}p_i\delta_i+(M+1).
$$
The first component in this sum is not smaller than
$$
(2m-M-1)(p_2+\dots+p_l)=(2m-M-1)p_1.
$$
For that reason the equality (\ref{7}) is impossible. Q.E.D. for
the proposition.\vspace{0.1cm}

The last step in the proof of Theorem 2 is the
following\vspace{0.1cm}

{\bf Proposition 7.} {\it The case $l>k$ is
impossible.}\vspace{0.1cm}

{\bf Proof.} Assume the converse: $l>k$. As in the proof of the
previous proposition, for any $i\leqslant k, i\geqslant 2$ we have
$$
B_{i-1}\subset S^{i-1}_1\cap E^{i-1}_1\cap E_{i-1}=(S^1\cap
E_1)^{i-1}\cap E_{i-1}.
$$
Set $\Delta=S^1\cap E_1$. Consider the hypersurface
$Q\subset{\mathbb P}$, containing the point $o$, which in the
affine coordinates $z_1,\dots,z_M$ is defined by the equation
$$
\widetilde{f}(z_*)=q_{m-1}(z_*)+\widetilde{q}_m(z_*)=0,
$$
where $q_{m-1}(z_*)$ is the same polynomial as in the equation
(\ref{8}) of the hypersurface $S$, and $\widetilde{q}_m(z_*)$ is a
generic homogeneous polynomial of degree $m$. Obviously, $Q^1\cap
E_1=\Delta$, the hypersurface $Q^1\subset X_1$ is non-singular and
the intersection of $Q^1$ with $E_1$ is everywhere transversal.
Therefore for every $i\leqslant k$, $i\geqslant 2$, we have
$$
B_{i-1}\subset\Delta^{i-1}\cap E_{i-1}=Q^{i-1}\cap E^{i-1}_1\cap
E_{i-1}.
$$
On the other hand, by the definition of the number $k$ we have
$B_k\not\subset S^k$, so that, because of the polynomial
$\widetilde{q}_m$ being general, we have $B_k\not\subset Q^k$.
Thus
$$
\mathop{\rm ord}\nolimits_T\varphi^*Q=\mathop{\rm
ord}\nolimits_{E_N}\varphi^*_{N,0}Q=(m-1)p_1+p_2+\dots+p_k=b.
$$
Now we argue ia the word for word the same way as in the proof of
Proposition 5: for the strict transform
$\widetilde{Q}\subset\widetilde{\mathbb P}$ we get the
presentation (\ref{9}), which immediately implies that
$(\widetilde{L}\cdot\widetilde{Q})=0$ for a general line $L\subset
{\mathbb P}$, which is impossible since the polynomial
$\widetilde{q}_m$ is a general one. Q.E.D. for the
proposition.\vspace{0.1cm}

Proof of Theorem 2 is complete.\vspace{0.3cm}


{\bf 6. Automorphisms of the hypersurface $S$.} Let us prove
Theorem 1. By Theorem 2 we need to show only the claim that the
group $\mathop{\rm Aut} ({\mathbb P}\setminus S)=\mathop{\rm Aut}
({\mathbb P})_S$ is finite, generically trivial. First of all,
every projective automorphism $\chi_{\mathbb P}$, preserving the
hypersurface $S$, maps the point $o$ to itself. Let $\mathop{\rm
Aut} ({\mathbb P})_o\subset \mathop{\rm Aut} ({\mathbb P})$ be the
stabilizer of the point $o$, and
$$
\pi\colon \mathop{\rm Aut} ({\mathbb P})_o\to \mathop{\rm Aut} (E)
$$
the natural projection, sending a projective automorphism $\xi\in
\mathop{\rm Aut} ({\mathbb P})_o$ to the corresponding
automorphism of the projectivized tangent space ${\mathbb
P}(T_o{\mathbb P})\cong E$. Obviously, for every $\chi_{\mathbb
P}\in \mathop{\rm Aut} ({\mathbb P})_S$ its image
$\pi(\chi_{\mathbb P})$ preserves the hypersurface $S^+\cap E$
(that is, the hypersurface $\{q_{m-1}=0\}$ in the sense of the
equation (\ref{8})). By \cite{MatMon}, the group $\pi(\mathop{\rm
Aut} ({\mathbb P})_S)$ is finite, and for a Zariski general
hypersurface $S$, trivial. Setting
$$
\pi_S=\pi|_{\mathop{\rm Aut} ({\mathbb P})_S},
$$
we see that it is sufficient to show that the kernel $\mathop{\rm
Ker}\pi_S$ is trivial. This is really easy.\vspace{0.1cm}

Every projective automorphism $\xi\in \mathop{\rm Ker}\pi$ in a
system of homogeneous coordinates $(x_0:x_1:\dots :x_M)$, such
that $o=(1:0:\dots :0)$, has the form
$$
\xi\colon (x_0:x_1:\dots :x_M)\mapsto (a_0x_0+a_1x_1+\dots
+a_Mx_M:x_1:\dots :x_M),
$$
where $a_0\neq 0$. The hypersurface $S$ in such a system of
coordinates is given by the equation $\Phi(x_*)=0$, where
$$
\Phi(x_0,\dots,x_M)=x_0q_{m-1}(x_1,\dots,x_M)+q_m(x_1,\dots,x_M),
$$
see the equation (\ref{8}). If $\xi\in \mathop{\rm Ker}\pi_S$,
then the homogeneous polynomial
$$
(a_0x_0+a_1x_1+\dots
+a_Mx_M)q_{m-1}(x_1,\dots,x_M)+q_m(x_1,\dots,x_M)
$$
is proportional to $\Phi(x_*)$. It is easy to see that this is
possible in one case only, when $a_0=1$ and $a_1=\dots =a_M=0$,
that is, $\xi=\mathop{\rm id}\nolimits_{\mathbb P}$. This
completes the proof of Theorem 1.


\begin{flushleft}
Department of Mathematical Sciences,\\
The University of Liverpool
\end{flushleft}

\noindent{\it pukh@liverpool.ac.uk}

\end{document}